\documentclass[%
 reprint,
 amsmath,amssymb,
 aps, 
pre,
]{revtex4-1}

\usepackage{float}
\usepackage{subfigure}
\usepackage{caption}
\usepackage{graphicx}
\usepackage{dcolumn}
\usepackage{bm}

\begin{document}

\preprint{APS/123-QED}

\title{Evolutionary game dynamics of controlled and automatic decision-making}

\author{Danielle F. P. Toupo}%
\email{dpt35@cornell.edu}
\affiliation{Center for Applied Mathematics, Cornell University, Ithaca, New York 14853\\
}%
\author{Steven H. Strogatz}%
\email{strogatz@cornell.edu}
\affiliation{Center for Applied Mathematics, Cornell University, Ithaca, New York 14853\\
}%
\author{Jonathan D. Cohen}%
\email{jdc@princeton.edu}
\affiliation{Department of Psychology, Princeton University, Princeton, New Jersey 08540\\
}%
\author{David G. Rand}%
\email{david.rand@yale.edu}
\affiliation{Department of Psychology, Yale University, New Haven, CT 06511}%

\date{\today}

\begin{abstract}
We integrate dual-process theories of human cognition with evolutionary game theory to study the evolution of automatic and controlled decision-making processes. We introduce a model where agents who make decisions using either automatic or controlled processing compete with each other for survival. Agents using automatic processing act quickly and so are more likely to acquire resources, but agents using controlled processing are better planners and so make more effective use of the resources they have. Using the replicator equation, we characterize the conditions under which automatic or controlled agents dominate, when coexistence is possible, and when bistability occurs. We then extend the replicator equation to consider feedback between the state of the population and the environment. Under conditions where having a greater proportion of controlled agents either enriches the environment or enhances the competitive advantage of automatic agents, we find that limit cycles can occur, leading to persistent  oscillations in the population dynamics. Critically, however, these limit cycles only emerge when feedback occurs on a sufficiently long time scale. Our results shed light on the connection between evolution and human cognition, and demonstrate necessary conditions for the rise and fall of rationality.
\end{abstract}

\maketitle

{\bf Dual-process theories of human cognition play a central role in the behavioral sciences. According to these theories, decisions are often made using either automatic processes which are fast and effortless but focused on the present, or controlled processes which are slow and effortful but can plan for the future. Evolutionary game theory models, however, almost never consider these distinctions. Therefore, little is known about the evolutionary dynamics of automatic versus controlled processing. Here, we address this gap by introducing an analytically tractable model for the evolution of agents that use automatic or controlled processing. The agents both compete with each other and alter their shared environment. We show that under certain circumstances, automatic and controlled processing can stably coexist within the population. We also identify conditions under which limit cycles occur. In such cases, the success of controlled agents alters the environment in a way that allows automatic agents to invade, and vice versa. Our results help to explain why human evolution may not necessarily be characterized by ever-increasing levels of rationality and forward-thinkingness, but instead may recurrently fall prey to periods of myopia.}

\section{Introduction}
Dual-process theories of human decision-making conceptualize decisions as arising from the interaction of (i) automatic processes that are ``hardwired" and thus computationally efficient but rigid; and (ii) controlled processes that are effortful but flexible~\cite{kahnemanbook,cohen1990control,miller2001integrative,posner1975attention,schneider1977controlled,shiffrin1997controlled,barett2006modularity}. Such a perspective has proved useful for understanding behavior across a wide range of domains, and has been used heavily in fields such as neuroscience~\cite{mcclure2004seperate,hare2009self}, cognitive and social psychology~\cite{evans2003two,evans2013dual,stanovich2000individual,tversky1983extensional,rand2012spontaneous,crockett2013models,cushman2013action}, and behavioral economics~\cite{fudenberg2006dual,kahneman2003maps,thaler1981economic}. Yet, despite playing a key role in human evolution, the interaction (and conflict) between automatic and controlled processing has been almost entirely overlooked by evolutionary game theorists.

Controlled processing is a defining feature of human cognition, thought to underlie virtually all higher level, characteristically human cognitive functions, such as planning, problem-solving reasoning, and symbolic language -- functions that, at least under some conditions, are capable of identifying and executing rational and even optimal behavior.  This might be taken to suggest that evolution should favor controlled processing, and that given sufficient time, control should prevail as the dominant mode of cognition.  However, there is evidence that human history is characterized by cyclical dynamics that suggest a proliferation of behaviors and social structures reflective of controlled processing, only to be followed by their demise and collapse~\cite{diamond2005collapse,richerson2009cultural}.  What might explain these historical cycles? Here, we explore the possibility that they may reflect the dynamics of interaction between automatic and controlled processing at the population level.  We do so by integrating dual-process agents into an evolutionary game-theoretic framework.

We focus our investigation of automatic versus controlled processing on a particular cognitive function: intertemporal choice~\cite{ainslie1975specious,laibson1997golden,thaler1981empirical}.  Intertemporal choice refers to decisions between options or behaviors that yield immediate reward versus those that are less rewarding in the short run, but have the potential to yield greater reward in the future.  

We choose to focus on intertemporal choice for three reasons.  First, the prevalence of short-sighted behavior has been identified as an important contributory factor to the demise of advanced civilizations~\cite{diamond2005collapse}, and is a topic of modern concern (e.g., failures to save, over-consumption of environmental resources, and abuse of antibiotics).  Second, immediacy-biased behaviors have been linked to automatic processing, while future-oriented behaviors have been linked to the engagement of controlled processing, both at the behavioral and neural levels of analysis~\cite{fudenberg2006dual,thaler1981economic,hare2009self,ward2000don}.  Thus, intertemporal choice may be a useful probe for studying the consequences of interactions between automatic and controlled processing at the population level.  Third, we have performed preliminary computer simulations that support this suggestion~\cite{rand2015}. In these simulations, agents foraged for resources (e.g., food) in an environment, and either consumed found resources immediately (when using automatic processing), or according to an optimal consumption plan calculated using a complex algorithm based on past experience.  Intriguingly, these simulations sometimes gave rise to evolutionary cycles in which the proportion of controlled agents in the population waxed and waned periodically. However, the complexity of the model led to analytical intractability, making it hard to understand what conditions gave rise to these cyclical dynamics, and what factors were responsible for the oscillations. 

A desire to understand these issues led us to the simplified model proposed in this paper. Using the replicator equation, a nonlinear dynamical system studied in evolutionary game theory~\cite{hofbauer1998evolutionary, nowak2006book}, we introduce a minimal model of dual-process agents engaged in intertemporal choice that captures the critical features of the scenario above while remaining sufficiently simple to be mathematically tractable.  In doing so, we provide a formal characterization of the conditions under which cyclical dynamics emerge, and the forces that drive such cycles.


\section{The model}
We model a world in which agents forage for goods, compete for access to these goods, and choose how to consume goods they acquire to generate fitness, with fitness being subject to diminishing marginal returns on consumption.  Agents are then subject to natural selection based on their resulting fitnesses. 

For simplicity, we assume there are only two types of agents, fully controlled and fully automatic, and we explore the evolution of the fraction of controlled agents, denoted as $x$. Automatic agents differ from controlled agents in two ways: how likely they are acquire goods (where the speed and efficiency of automaticity is advantageous), and how they choose to consume those resources (where the rationality and planning ability of control is advantageous). 

The world is parametrized by the probability $\rho$ of finding a good (all goods are of equal size, normalized to $1$ energy unit), and the competitive advantage $\beta$ that automatic agents have over controlled agents in acquiring goods (where $\beta = 0$ means that both types of agents have an equal probability of acquiring goods).

\subsection{Competitive advantage}
Because automatic processing is assumed to be faster and less taxing than controlled processing, automatic agents have a competitive advantage over controlled agents when seeking to acquire goods. For example, it could be that when agents of both types simultaneously encounter a good, the automatic agent acts more quickly and snatches the good before the controlled agent can respond. Or it could be that the ponderous deliberation engaged in by controlled agents sometimes causes them to miss opportunities that an automatic agent would be more likely to exploit. 

As a result, automatic agents are more likely to acquire a good in any given time period, and so the two types of agents differ in their expected waiting time between acquiring goods (i.e., the average number of time steps between acquiring one good and the next). We define the probability of acquiring a good as $p_A$ for automatic agents, and $p_C$ for controlled agents. Thus the average waiting time for an automatic agent $\tau_A$ is given by\\
\begin{equation}
\tau_A = \frac{1}{p_A},\;\; \text{with} \;\; p_A = \rho (1+ \beta x),
\label{eqn:pA}
\end{equation}
while the average waiting time for a controlled agent $\tau_C$ is given by
\begin{equation}
\tau_C = \frac{1}{p_C},\;\; \text{with}\;\; p_C = \rho(1- \beta (1-x)).
\label{eqn:pC}
\end{equation}

For $\rho > 0$ and $\beta > 0$, it is the case that $\tau_A < \tau_C$: automatic agents acquire goods more frequently than controlled agents (again, because automatic processing is faster and more efficient). 

\subsection{Consumption}
To implement diminishing marginal returns on resource consumption, we define the fitness gained from consuming a fraction $z$ of a good as $z/(a+z)$, where $a$ controls the extent of diminishing marginal returns, with lower $a$ leading to more steeply diminishing returns. Recall that goods are normalized to have size $1$ when acquired. 

When automatic agents acquire a good, they consume all of it immediately; hence $z=1$, yielding a fitness benefit of $1/(a+1)$. They then spend, on average, the next $\tau_A - 1$ time steps consuming nothing, until they again acquire a good. Therefore the expected fitness per time step of an automatic agent is given by
\begin{equation}
\displaystyle f_A = \frac{\frac{1}{1+a}}{\tau_A} =  \frac{\rho +\beta  \rho  x}{a+1},
\label{eqn:fA}
\end{equation} 
from Eq.~\eqref{eqn:pA}.

In contrast, controlled agents consume acquired resources more carefully: they pace their consumption, spreading it out evenly so as to obtain the maximum possible amount of fitness gain from it. (Because of the diminishing marginal returns on consumption, it is wasteful to consume the entire resource immediately; evenly spaced consumption results in greater fitness yield.) Thus, the prudent planning of controlled agents leads them to consume $z = 1/\tau_C$ units of good in each of the $\tau_C$ time steps, and thereby to gain a fitness benefit 
\begin{equation}
f_C =\frac{\frac{1}{\tau_C}}{a+\frac{1}{\tau_C}} = \frac{\rho  (\beta  (x-1)+1)}{a+\rho  (\beta  (x-1)+1)},
\label{eqn:fC}
\end{equation}
from Eq.~\eqref{eqn:pC}.

\section{Evolutionary dynamics in a constant environment}
Having defined the fitness of the two types of agents, we turn to evolutionary dynamics. Specifically, we ask which strategy (or combination of the two) will be favored by natural selection for different fixed values of resource availability $\rho$ and competitive advantage of automatic agents $\beta$. We do so using the replicator equation from evolutionary game theory~\cite{hofbauer1998evolutionary, nowak2006book} to characterize how the relative fractions of controlled and automatic agents, $x$ and $1-x$, respectively, vary over time. The replicator equation compares the fitness of controlled agents to the population average fitness. It increases the frequency of controlled agents over time if they have higher fitness than automatic agents, and decreases it if the opposite is true. 

The replicator equation for our system, using \eqref{eqn:fA} and \eqref{eqn:fC},  is given by
\begin{eqnarray}
\dot{x} &=& x(f_C-(x f_C +(1-x) f_A)) \nonumber \\
             &=&(x-1) x \left(\frac{a}{a-\beta  \rho +\rho +\beta  \rho  x}+\frac{\rho +\beta  \rho  x}{a+1}-1\right).  \label{eqn:xdot}
\end{eqnarray} 
Note that we do not need a separate equation for the fraction of automatic agents because that quantity is given by $1-x$. The long-term dynamics of \eqref{eqn:xdot} are characterized in Fig.~1(a), where for the sake of illustration we fix $a= 0.15$ and vary $\beta$ and $\rho$. We see that the $(\beta, \rho)$ space is subdivided into five distinct regions. We describe the dynamics within each region below.

\begin{figure}[h!]
\center
\includegraphics[scale = 0.6,width=8.5cm]{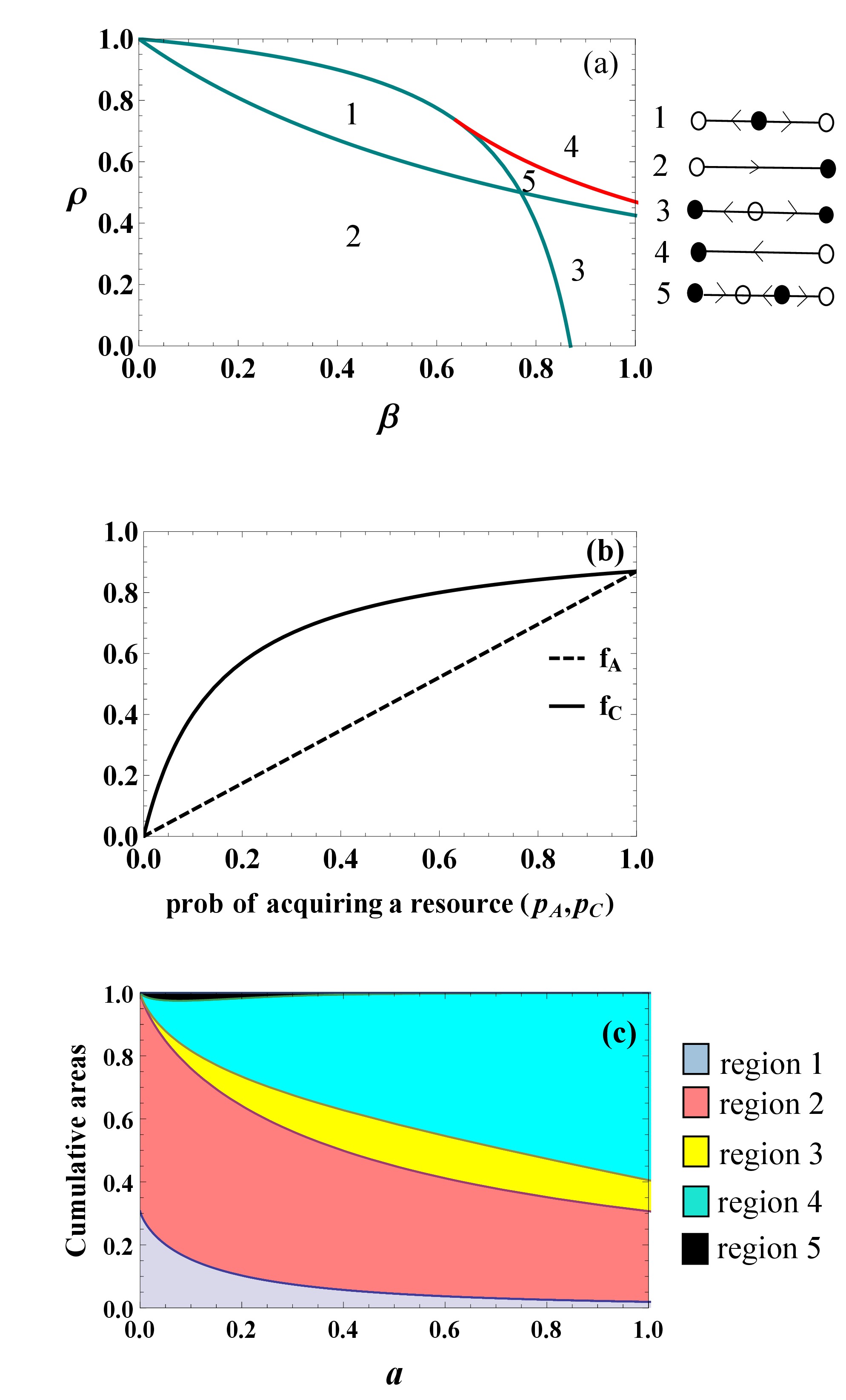}
\caption{\label{fig:fig21} (Color online) Bifurcation analysis of Eq.~\eqref{eqn:xdot}. (a) Stability diagram (left) and phase portraits (right) for Eq.~\eqref{eqn:xdot} with $a = 0.15$. Transcritical bifurcation, green curves; Saddle-node bifurcation, red curve. (b) Fitnesses $f_C$ and $f_A$ as functions of $p_A$ and $p_C$, for $a = 0.15$. (c) Areas of regions (1)-(5) in the stability diagram, as function of $a$.}
\end{figure}

The endpoint solutions of $x = 0$ (all automatic agents) and $x = 1$ (all controlled agents) are always fixed points regardless of $\beta$ and $\rho$. In regions $2$ and $4$, these are the only fixed points. When resources are scarce and the competitive advantage of automatics is low (region $2$), $x=1$ is the global attractor and control dominates automatic processing. Conversely, when resources are plentiful and the competitive advantage of automatics is high (region $4$), $x=0$ is the global attractor and automatic processing dominates control. This is because on the one hand, automatic agents always consume $\displaystyle \rho \beta$ more goods on average than controlled agents in each time step (given \eqref{eqn:pA} and \eqref{eqn:pC}); but on the other hand, controlled agents make more judicious use of those resources (as controlled by \textit{a}). Therefore, for a given value of \textit{a}, control wins when $\rho \beta$ is sufficiently small and automaticity wins when $\displaystyle  \rho \beta$ is large. The smaller \textit{a} is (i.e., the greater the diminishing marginal returns on consumption), the larger region $2$ is and the smaller region $4$ is. 

In the other regions, however, there can be up to two interior fixed points, in addition to these endpoint solutions. The first results from having a relatively resource-rich world with relatively little competitive advantage of automatics (regions 2 and 5). This interior fixed point is always stable, and leads to coexistence of automatic and controlled processing. The second results from a relatively resource poor world in which the competitive advantage of automatics is relatively common (regions 3 and 5). This interior fixed point, by contrast, is always unstable and leads to bistability between automatic and controlled processing. 

To understand why a rich world with little competitive advantage for automatics leads to coexistence while a poor world with high competitive advantage for automatics leads to bistability, we must consider how selection pressure varies based on the makeup of the population. In general, coexistence occurs when each strategy is at an advantage when it is rare, whereas bistability occurs when each strategy is at a disadvantage when it is rare. From \eqref{eqn:pA} and \eqref{eqn:pC}, we see that increasing the fraction of controlled agents $x$ by a given amount also increases the probability of finding a resource for both types of agents ($p_A$ and $p_C$) equally, regardless of $\rho$ and $\beta$. Therefore, what determines the dynamics when $x$ is small versus large is how $x$ (and the resulting increase in the probability of finding a good) translates into fitness for automatic versus controlled agents (which does depend on $\rho$ and $\beta$). From  \eqref{eqn:fA} and \eqref{eqn:fC}, we see that $f_A$ is linear in $p_A$, whereas $f_C$ is a nonlinear function of $p_C$ (see Fig.1(b)). Thus, because of the concavity of $f_C$, an increase in the fraction of controlled agents can have different effects on the relative fitness of automatic versus controlled processing depending on $\rho$ and $\beta$.

In a rich world ($\rho$ large) with relatively weak competitive advantage for automatics ($\beta$ small), as found in region 1, resources are common and $p_C$ and $p_A$ are relatively close to $1$. Thus the dynamics sit in a region where the $f_C$ curve in fig. 1.b has a shallower slope than that of the linear $f_A$ curve. Consequently, going from $x=0$ to $x=1$ leads to a bigger increase in fitness for automatic agents than controlled agents. As a result, this produces a situation in which (with the right $\rho$ and $\beta$) control outperforms automatic near $x=0$ (when control is rare), but as $x$ increases, the advantage of control dissipates and reverses such that automatic outperforms control near $x=1$. Thus neither endpoint is stable, leading to coexistence.

Different dynamics occur in a poor world ($\rho$ small) where the competitive advantage for automatics is high and $\beta$ is large (region 3).  Here, $p_C$ and $p_A$ are relatively close to $0$. In this case, the slope of $f_C$ is larger than that of $f_A$, and thus going from $x=0$ to $x=1$ leads to a greater increase in fitness for controlled agents than automatic agents. This produces a situation in which automatic agents outperform controlled agents near $x=0$, whereas controlled agents outperform automatic agents near $x=1$. Thus both endpoints are stable, leading to bistability. 

Finally, when both $\rho$ and $\beta$ are moderately high (region 5), the resulting long-term dynamics are a mix of regions 1 and 2, with bistability occurring between $x=0$ and a stable interior fixed point (i.e., coexistence).

Fig.~1(c)  shows the areas of the five regions in Fig.~1(a) as the parameter $a$ increases. Increasing $a$ increases the size of the region in $(\beta, \rho)$ space where automatic agents dominate (region 4) and drastically decreases the regions corresponding to bistability, coexistence, or dominance of controlled agents (regions 1, 2, 3, and 5).

\section{Feedback between the population and the environment}
\noindent In Section III we assumed that the environment is constant, such that the parameters $\rho$ and $\beta$ are fixed. There are many situations, however, in which the current makeup of the population can influence the environment, often with some lag~\cite{cohen2005vulcanization,diamond2005collapse,richerson2009cultural}. Thus in this section, we extend the model from Section III to incorporate such feedback effects. To do so, we introduce a modified version of the replicator equation that includes additional differential equations describing how the environmental parameters $\beta$ and $\rho$ vary with $x$, the fraction of controlled agents in the population. In Section IV.A, we analyze a system in which an increase in controlled processing increases $\beta$, thus augmenting the competitive advantage of automatic agents.  In Section IV.B, we analyze a system in which an increase in controlled processing increases $\rho$. This scenario models a situation where greater use of controlled processing enriches the environment and enhances resource availability for everyone, thanks (for example) to increased technological innovation leading to greater agricultural output. In Section IV.C, we analyze a system with both of these features.

\subsection{Scenario 1: Controlled processing increases competitive advantage of automaticity}
\noindent Here we consider the consequences of allowing $\beta$ to positively co-vary with $x$. This implements a scenario in which having more controlled agents leads to greater population density and thus a larger $\beta$. The increase in population density could reflect larger population size, which results directly from the fact that populations with more controlled agents have higher average fitness. Alternatively, it could reflect an externality such as cognitive control allowing people to live more densely without violent conflict. 

To link $\beta$ and $x$, we introduce a differential equation on $\beta$ that pulls its value towards the current value of $x$. We also incorporate the possibility of lag, specified by a parameter $\tau_{\beta}$. This lag captures the effect that an increase in $x$ at time $t$ does not always have an immediate impact on $\beta$. For example, increased birth rates do not immediately lead to larger numbers of competing adults. 

The new system is given by
\begin{eqnarray}
\dot{x} &=& x(f_C-(x f_C +(1-x) f_A)) \nonumber \\ 
\dot{\beta} &=& \frac{x-\beta}{\tau_{\beta}}. \label{eqn:xbeta1}
\end{eqnarray}
After insertion of \eqref{eqn:fA} and \eqref{eqn:fC}, \eqref{eqn:xbeta1} becomes
\begin{eqnarray}
\dot{x} &=& (x-1) x \left(\frac{a}{a-\beta  \rho +\rho +\beta  \rho  x}+\frac{\rho +\beta  \rho  x}{a+1}-1\right). \nonumber  \\ 
\dot{\beta} &=& \frac{x-\beta}{\tau_{\beta}} \label{eqn:xbeta2}
\end{eqnarray}

Note that the $\dot{x}$ replicator equation is the same as it was previously, except that now $\beta$ is also a variable. Also note that in equilibrium, $\beta=x$.         

To illustrate how this addition of feedback between the population and the environment affects the dynamics, we begin by fixing $a=0.8$ and examining the effect of $\rho$ and $\tau_{\beta}$ (Fig. 2(a) and 3). We find three possible types of long-term dynamics: dominance of controlled agents (region 1, no interior fixed point); coexistence of automatic and controlled agents (region 3, stable interior fixed point); and limit cycles in which both types of agents are present but their relative abundances oscillate (region 2, unstable interior fixed point). 

The parameter regime in Fig.~2(a) for which limit cycles exist is bounded by two vertical asymptotes, and within that strip, $\tau_{\beta}$ must be sufficiently large. In Fig.~2(a), $a = 0.8$ and limit cycles exist if $0.1 < \rho < 0.52$ and $\tau_{\beta} > 104.47$. Specifically, the limit cycles are born in a supercritical Hopf bifurcation. The equation of the Hopf bifurcation curve has been computed analytically and is too complicated to show.

\begin{figure}[h!]
\centering
\includegraphics[scale=0.7,width=6.8cm]{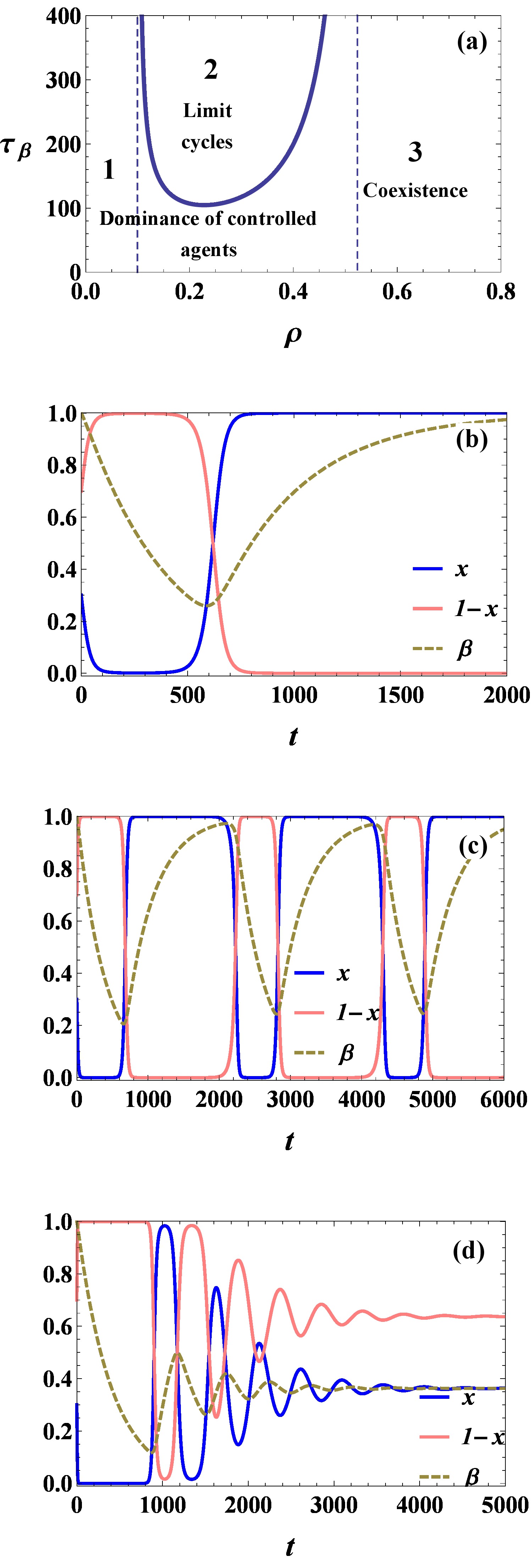}
\label{fig:fig2}

\caption{(Color online) Bifurcation analysis of Eq.~\eqref{eqn:xbeta2} with $a=0.8$ and $\tau_{\beta} = 400$. (a) Stability diagram. Hopf bifurcation, blue curve. (b) Time series of a typical solution in region 1 with $\rho = 0.1.$ (c) Time series of a typical solution in region 2 with $\rho = 0.2.$ (d) Time series of region 3 with $\rho = 0.65$.}
\end{figure}

\begin{figure}
\centering
\includegraphics[scale = 0.6,width=7cm]{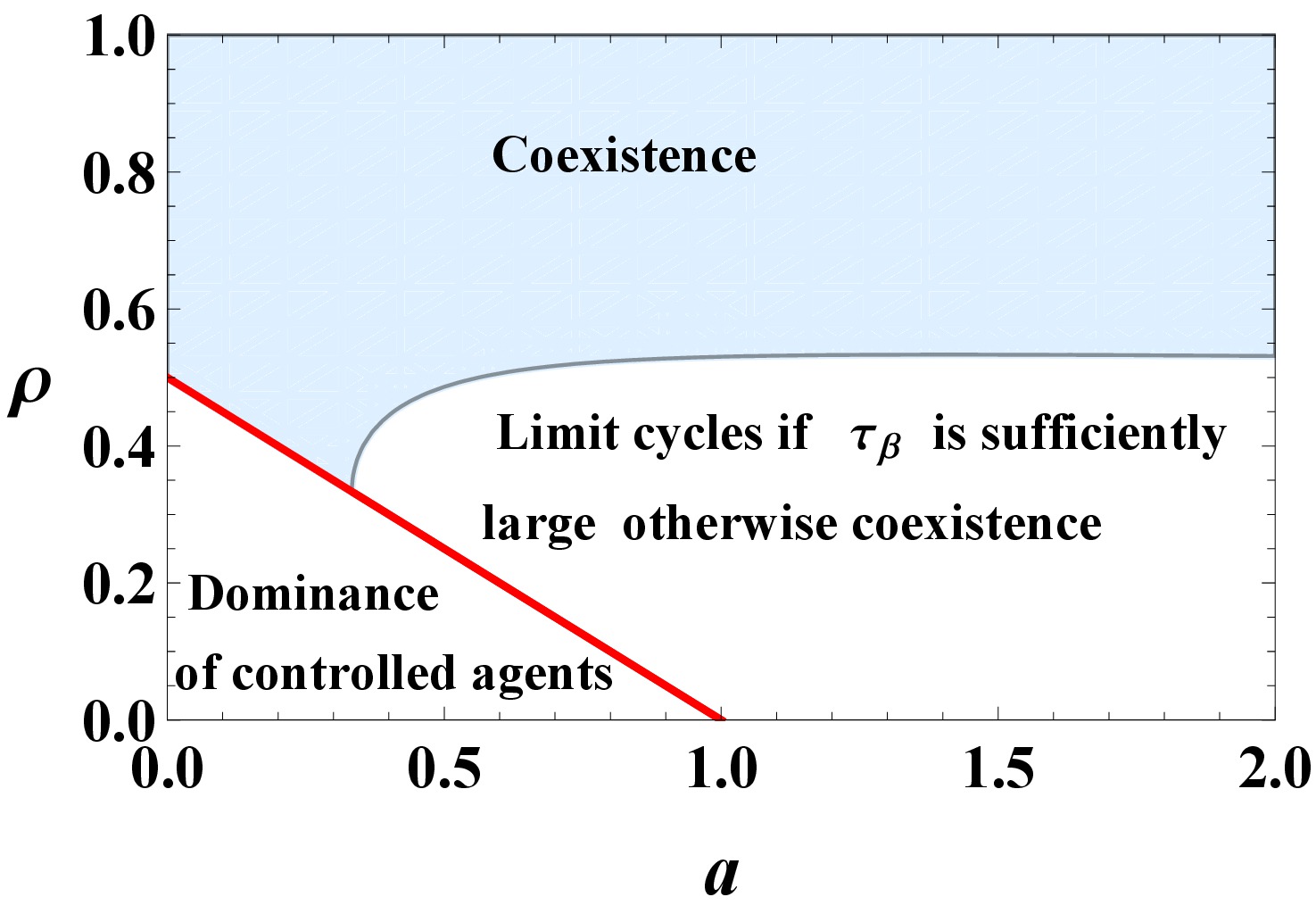}

\caption{(Color online) Characterization of the $(x,\beta)$ system from Eq.~\eqref{eqn:xbeta2}}
\end{figure}

We conclude this section by asking how $a$, the extent of diminishing marginal returns on consumption, changes the dynamics of Eq.~\eqref{eqn:xdot} (Fig. 3). We find that only the three types of dynamics observed in Figs.~3 and Figs.~2(b)-(d) are possible: if $\rho < (a+1)/2$, the long-term behavior is dominance of controlled agents (no interior fixed point); if $\rho > (a+1)/2$ and $\rho > \rho^*$, the long-term behavior will be coexistence; and if $\rho > (a+1)/2$ and $\rho < \rho^*$, limit cycles are possible if $\tau_{\beta}$ is sufficiently large, otherwise there will be coexistence. The curve bounding the region in $(a, \rho)$ space where limit cycles are possible has been computed numerically.  

In sum, we see that limit cycles can arise from feedback between the overall population density and the fraction of controlled agents in the population. Critically, these oscillations emerge only when the feedback is sufficiently delayed, in which case they occur over a wide range of $\rho$ and $a$ values. 

\subsection{Scenario 2: Controlled processing increases resource availability}
\noindent Here, we leave $\beta$ fixed and instead link $\rho$ to $x$, using the same formulation for $\rho$ here as for $\beta$ in Scenario 1. This models a scenario in which controlled agents enrich the environment, say by creating technologies that increase resource abundance for everybody. Again, we add a lag that represents the time required for the development of such technologies and their impact on the environment to occur. This gives rise to the following system:

\begin{eqnarray}
\dot{x} &=& x(f_C-(x f_C +(1-x) f_A)) \nonumber \\ 
\dot{\rho} &=& \frac{x-\rho}{\tau_{\rho}}. \label{eqn:xrho1}
\end{eqnarray}
After insertion of \eqref{eqn:fA} and \eqref{eqn:fC}, Eq.~\eqref{eqn:xrho1} becomes
\begin{eqnarray}
\dot{x} &=& (x-1) x \left(\frac{a}{a-\beta  \rho +\rho +\beta  \rho  x}+\frac{\rho +\beta  \rho  x}{a+1}-1\right) \nonumber  \\ 
\dot{\rho} &=& \frac{x-\rho}{\tau_{\rho}}. \label{eqn:xrho2}
\end{eqnarray}

Again, $\dot{x}$ is the same in the system without feedback, and in equilibrium $\rho=x$. For the sake of illustration, we fix $a=1.5$ and examine the dynamics as a  function of $\beta$ and $\tau_{\rho}$ (Fig. 4(a)). We find three possible types of long-term dynamics: coexistence of automatic and controlled agents (region 1, stable interior fixed point); limit cycles in which both types of agents are present but their relative abundances oscillate (region 2, unstable interior fixed point), and dominance of automatic agents (region 3, no interior fixed point). As Fig. 4(a) indicates, limit cycles exist only if $0.249 < \rho < 0.4$ and $\tau_{\rho} > 446.3$. The Hopf bifurcation curve bounding the limit cycle region has been calculated analytically, but is too complicated to show here. Figures 4(b)-(d) show time series for sample trajectories from within each region. 

\begin{figure}[h!]
\centering
\includegraphics[width = 6.9cm]{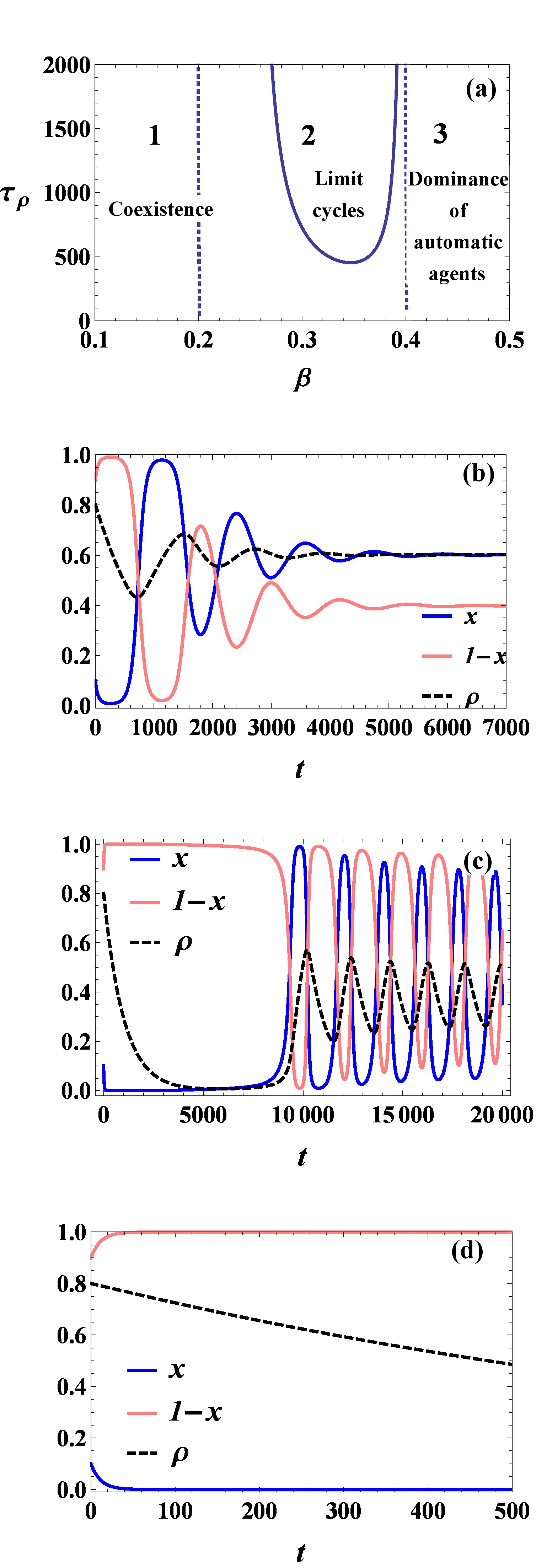}
\label{fig:fig4}

\caption{(Color online) Bifurcation analysis of Eq.~\eqref{eqn:xrho2} with $a=1.5$ and $\tau_{\rho} = 1000$. (a) Stability diagram. Hopf bifurcation, blue curve. (b) Time series of region 1 with $\beta = 0.2.$ (c) Time series of region 2 with $\beta = 0.3.$ (d) Time series of region 3 with $\beta = 0.45$.}
\end{figure}

\begin{figure}
\centering
\includegraphics[width=6cm]{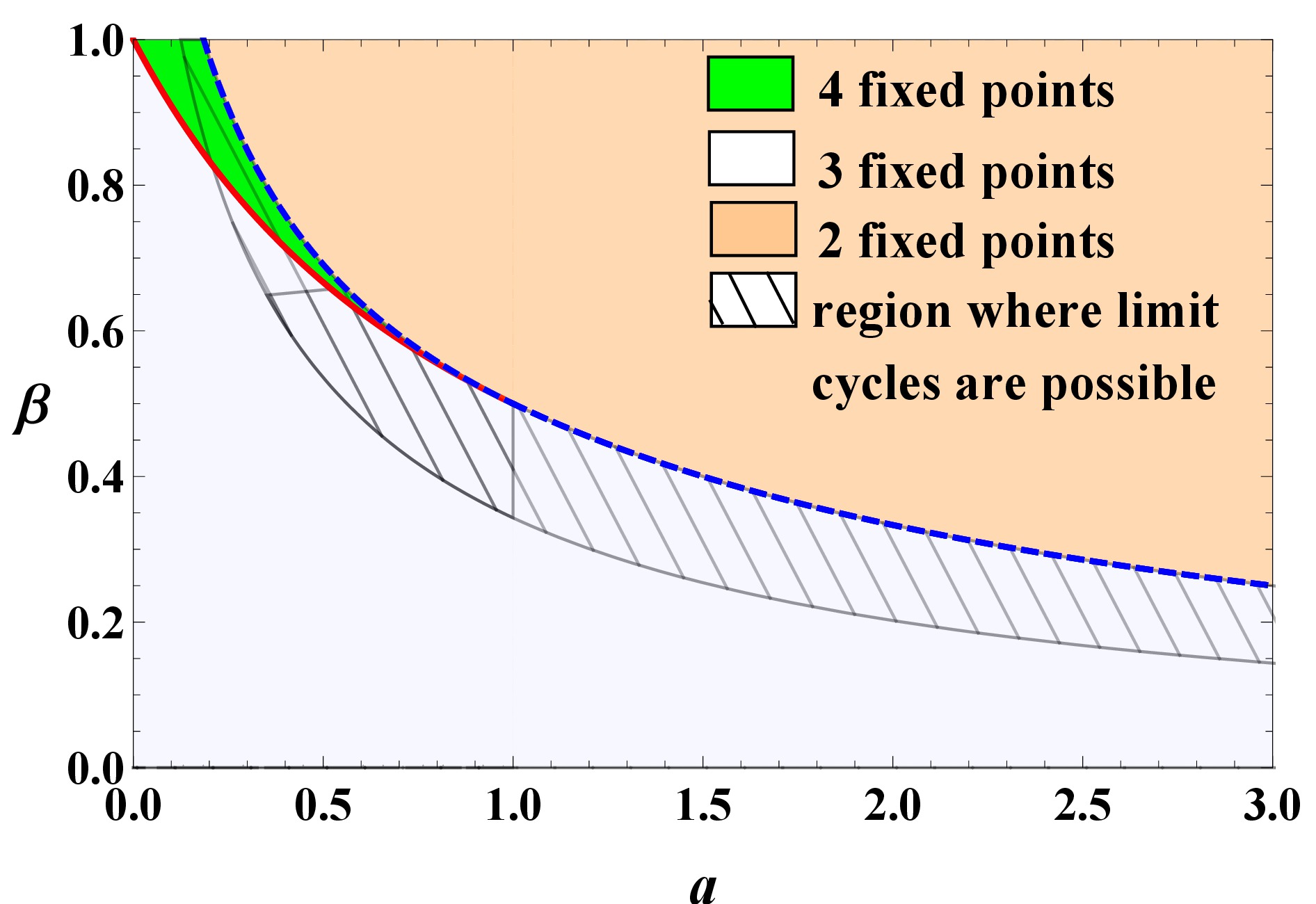}

 \caption{(Color online) Characterization of \eqref{eqn:xrho2} with $a=1.5.$}
\end{figure}

\begin{figure}
\centering
\includegraphics[width=6cm]{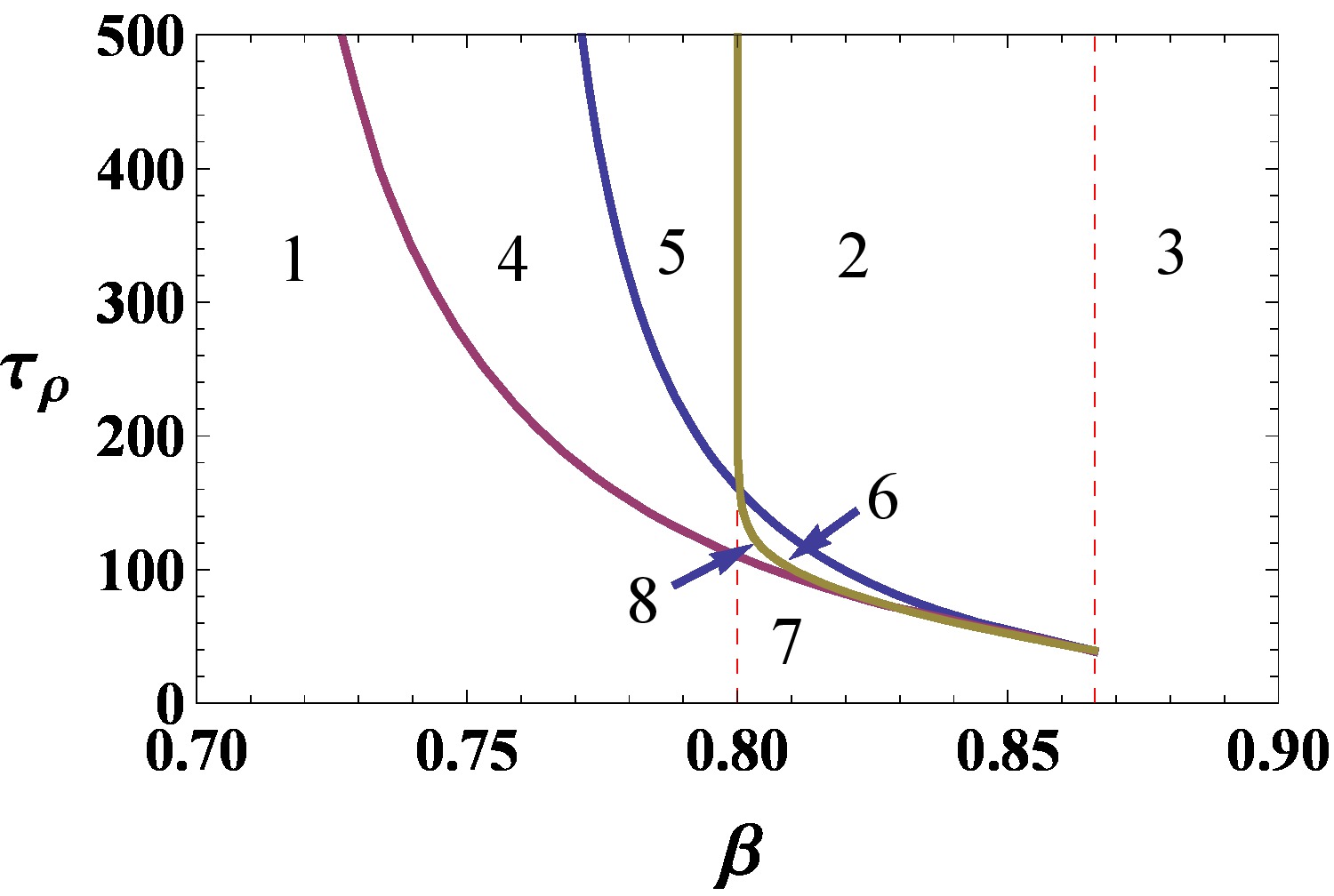}

\caption{(Color online) Bifurcation diagram of \eqref{eqn:xrho2} with $a = 0.5$. Hopf bifurcation, blue curve; fold bifurcation (saddle-node coalescence) of cycles, purple curve; Homoclinic bifurcation, orange curve.}
\end{figure}

There is an important difference between these dynamics and the dynamics studied in Scenario 1 when $x$ and $\beta$ were positively correlated: in Scenario 1, dominance of controlled but not automatic agents was possible, whereas here in Scenario 2, the opposite is true.  Only automatic agents can dominate.  Moreover, this dominance by automatic agents occurs only if $a > 1$.

Next we ask how the dynamics of ~\eqref{eqn:xrho2} depend on the parameter $a$, which reflects the strength of diminishing returns. We find that only the three types of long-term dynamics observed in Fig. 4(a) are possible if $a > 1$, but that more complex dynamics emerge when $a < 1$ (Figs. 5 and 6).  Figure 6 characterizes the dynamics as a function of $\beta$ and $\tau_\rho$ for $a = 1/2$. The long-term dynamics of ~\eqref{eqn:xrho2} for $a< 1$ sometimes depend on the initial conditions, in a manner that can be summarized as follows:
\begin{itemize}
\itemsep0em 
\item Regions 3 and 5: dominance of automatic agents.
\item Region 2: either limit cycle oscillations of the two strategies, or dominance of automatic agents, depending on the initial conditions.
\item Region 4:  either oscillations or coexistence, depending on the initial conditions.
\item Regions 6 and 7: either dominance of automatic agents or coexistence,  depending on the initial conditions.
\item Region 8:  oscillation of the two strategies, dominance of automatic agents, or coexistence, depending on the initial conditions. 
\end{itemize}

In summary, adding feedback between the fraction of the population that uses controlled processing and the availability of resources can also give rise to limit cycles when the feedback is sufficiently delayed. Compared to the $(x,\beta)$ system discussed in Scenario 1 (Section IV.A) , however, limit cycles occur over a smaller range of $(\beta,a)$ combinations. Furthermore, the dynamics of the $(x, \rho)$ system of Scenario 2 are substantially more complex than those of the $(x,\beta)$ system of Scenario 1.

\subsection{Scenario 3: Controlled processing increases both competition and resource availability}
\noindent Finally, we consider the case in which $x$ influences both $\beta$ and $\rho$. To do so, we use a three differential equation system that includes the original replicator equation for $\dot{x}$, as well as the $\dot{\beta}$ equation from ~\eqref{eqn:xbeta1} and the $\dot{\rho}$ equation from ~\eqref{eqn:xrho2}, with the use of two different time-constants, $\tau_{\rho}$ and $\tau_{\beta}$, for the two different environmental feedback equations. Thus our system is given by

\begin{eqnarray}
\dot{x} &=& x(f_C-(x f_C +(1-x) )) \nonumber \\ 
\dot{\beta} &=& \frac{x-\beta}{\tau_{\beta}} \nonumber \\
\dot{\rho} &=& \frac{x-\rho}{\tau_{\rho}}. \label{eqn:xetarho1}
\end{eqnarray}

After insertion of \eqref{eqn:fA} and \eqref{eqn:fC}, Eq.~\eqref{eqn:xetarho1} becomes

\begin{eqnarray}
\dot{x} &=&   (x-1) x \left(\frac{a}{a-\beta  \rho +\rho +\beta  \rho  x}+\frac{\rho +\beta  \rho  x}{a+1}-1\right)\nonumber  \\
\dot{\beta} &=& \frac{x-\beta}{\tau_{\beta}} \nonumber \\
\dot{\rho} &=& \frac{x-\rho}{\tau_{\rho}}. \label{eqn:xetarho2}
\end{eqnarray}

We perform numerical simulations to demonstrate that limit cycles can also arise in the 3D system, as shown in Fig. 7 with $a = 1.5$. We see that limit cycles are possible as long as neither $\tau_{\beta}$ nor $\tau_{\rho}$ are too small.

\begin{figure}[h]
\begin{centering}
\includegraphics[width=6cm]{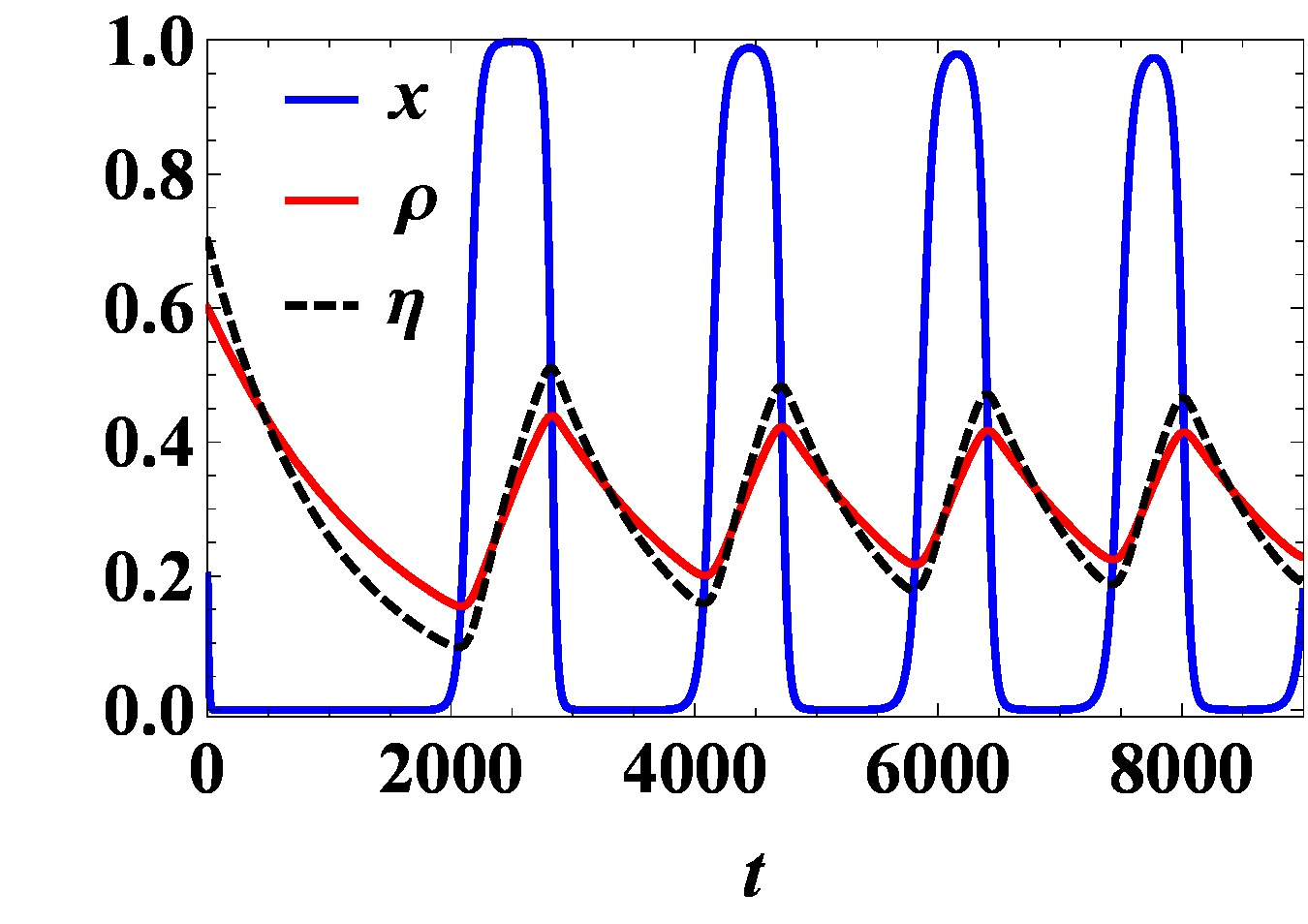}
\label{fig:fig7}

\caption{(Color online) Analysis of~\eqref{eqn:xetarho2} with $a = 1.5$,  $\tau_{\rho}= 1500$, $\tau_{\beta} = 1000$. (a) Parametric plot. (b) Time series.}
\end{centering}
\end{figure}

\section{Discussion}

Here we have introduced an analytically tractable model of the evolution of dual-process agents. Our model focuses on intertemporal choice, with agents foraging for, and competing over, goods which they consume to generate fitness. Agents that use automatic processing are at an advantage when acquiring goods because of their speed and efficiency, but immediately consume any goods they acquire in short-sighted fashion. Controlled agents, conversely, engage in long-term planning and make better use of the goods they manage to acquire. 

Within this framework, the agents' world is parametrized by $\rho$, the availability of resources (defined as the average probability of finding a good per unit time), and $\beta$, the competitive advantage of automatic agents (the increased likelihood of automatic agents acquiring a good over controlled agents). Our analysis allows us to characterize which parts of the ($\rho$, $\beta$) parameter space lead to dominance of automatic or controlled processing, bistability, or coexistence, as well as the conditions under which limit cycles arise.

In particular, we find that natural selection favors controlled agents when $\rho$ and $\beta$ are both small (poor worlds with little competition), automatic agents when $\rho$ and $\beta$ are both large (rich worlds with substantial competition), coexistence when $\rho$ is large and $\beta$ small (rich worlds with little competition), and bistability when $\rho$ is small and $\beta$ large (poor worlds with substantial competition). Furthermore, we find that limit cycles are a robust feature of adding environmental feedback whereby a greater frequency of controlled agents leads to either higher $\beta$, higher $\rho$, or both. Critically, however, the feedback must be sufficiently lagged in order for limit cycles to emerge.

Thus our analyses demonstrate the key role that feedback between the population and the environment plays in population (and ecological) dynamics. Such feedback can lead to cyclical dynamics that are otherwise impossible in a two-species competition model. Critically, environmental feedback is absent from typical evolutionary game-theoretic models, in which the game parameters are fixed, and only the population make-up varies over time~\cite{hofbauer1998evolutionary, nowak2006book}. By extending the replicator equation to include linkage between the population and one or more of the game parameters, we allow a richer range of dynamics that help to explain cyclical dynamics observed in human history. 

In the interest of analytical tractability, our model makes a number of simplifying assumptions. Most importantly, we consider the limiting case of entirely automatic agents competing with entirely controlled agents. In reality, agents exist on a continuum of inclination towards automaticity versus control.  We also consider a highly simplified foraging environment, and a simple decision rule for controlled agents (spread consumption out evenly over the expected waiting period until the next good is acquired). We are confident, however, that these particular simplifications did not distort our results, however, based on our prior computer simulation work~\cite{rand2015}. These simulations had agents that could engage in both automatic and controlled processing, and examined a much more complex foraging environment. Nonetheless, our simplified model recreates the same dynamics as the more complex simulations. 

The framework we introduce here can be extended in many ways to assess the impact of other simplifications, and to explore other questions. For example, spatial structure could be added~\cite{durrett1994importance,nowak2010evolutionary,perc2013}, agents could differ in the extent to which they impact the game parameters, or the game parameters could vary cyclically over time (instead of, or in addition to, variation caused by the population) ~\cite{rand2011evolutionary,ruelas2012nonlinear}. Our basic framework could also be applied to study dual-process cognition in domains beyond intertemporal choice, such as risky choice~\cite{zur1981effect,fudenberg2006dual} or cooperation in social dilemmas~\cite{rand2012spontaneous,rand2013human,rand2014social}.
In summary, we have introduced an evolutionary game-theoretic model of dual-process agents who make decisions using either automatic or controlled cognitive processing, and who not only compete with each other but also affect their environment. Our model demonstrates how the tendency for controlled processing to enrich the environment or grow the population undermines the advantages of controlled cognition, leading to the eventual invasion of automaticity and short-sightedness. Thus our model sheds light on historical cycles through which controlled processing, and associated phenomena such as careful planning and technological innovation, may rise and fall. The success of controlled cognition naturally leads to its own demise.

\section*{Acknowledgments}
This research was supported in part by a Sloan-Colman fellowship to DFPT, the National Science Foundation through grant DMS-1513179 to SHS, and the John Templeton Foundation through grants to JDC and DGR.  The authors also thank Richard Rand, Damon Tomlin, and Elliot Ludwig for helpful comments.

%
%

\end{document}